\documentclass[11pt]{article}

\usepackage{latexsym,epsfig}
\usepackage{amssymb,amsmath,amsfonts}
\usepackage{exscale} 
\usepackage{ifthen}
\usepackage{longtable}
\usepackage{subfigure}
\usepackage{caption}
\usepackage{cite}
\usepackage{lscape}

\begin{document}
\bibliographystyle{gabialpha}
\protect\pagenumbering{arabic}
\setcounter{page}{1}
 
\newcommand{\Zt}{\rm}

\newcommand{\ba}{\begin{array}}
\newcommand{\ea}{\end{array}}
\newcommand{\pot}{{\cal P}}
\newcommand{\curv}{\cal C}
\newcommand{\ddt} {\mbox{$\frac{\partial  }{\partial t}$}}
\newcommand{\hl}{\sf}
\newcommand{\hd}{\sf}

\renewcommand{\d}{\mathrm{d}}
\newcommand{\e}{\mathrm{e}}
\newcommand{\Ad}{\mbox{\rm Ad}}
\newcommand{\Adsm}{\mbox{{\rm \scriptsize Ad}}}
\newcommand{\ad}{\mbox{\rm ad}}
\newcommand{\adsm}{\mbox{{\rm \scriptsize ad}}}
\newcommand{\diag}{\mbox{\rm Diag}}
\newcommand{\sect}{\mbox{\rm sec}}
\newcommand{\id}{\mbox{\rm id}}
\newcommand{\idsm}{\mbox{{\rm \scriptsize id}}}
\newcommand{\eps}{\varepsilon}
\newcommand{\Summ}{P}
\newcommand{\ein}{\,\rule[-5pt]{0.4pt}{12pt}\,{}}

\newcommand{\aL}{\mathfrak{a}}
\newcommand{\bL}{\mathfrak{b}}
\newcommand{\mL}{\mathfrak{m}}
\newcommand{\kL}{\mathfrak{k}}
\newcommand{\gL}{\mathfrak{g}}
\newcommand{\nL}{\mathfrak{n}}
\newcommand{\hL}{\mathfrak{h}}
\newcommand{\pL}{\mathfrak{p}}
\newcommand{\uL}{\mathfrak{u}}
\newcommand{\lL}{\mathfrak{l}}

\newcommand{\kG}{{\tt k}}
\newcommand{\nG}{{\tt n}}

\newcommand{\Cart}{$G=K e^{\overline{\aL^+}} K$}
\newcommand{\Area}{\mbox{Area}}
\newcommand{\Hd}{\mbox{\rm Hd}}
\newcommand{\Hdim}{\mbox{\rm dim}_{\mbox{\rm \scriptsize Hd}}}
\newcommand{\Tr}{\mbox{\rm Tr}}
\newcommand{\bs}{{\cal \beta}}
\newcommand{\bv}{{\cal B}}

\newcommand{\nc}{{\cal N}}
\newcommand{\MM}{{\cal M}}
\newcommand{\Ch}{{\cal C}}
\newcommand{\clCh}{\overline{\cal C}}
\newcommand{\Sh}{\mbox{Sh}}
\newcommand{\smSh}{\mbox{{\rm \scriptsize Sh}}}
\newcommand{\Cnt}{\mbox{\rm C}}

\newcommand{\NN}{\mathbb{N}} \newcommand{\ZZ}{\mathbb{Z}}
\newcommand{\QQ}{\mathbb{Q}} \newcommand{\RR}{\mathbb{R}}
\newcommand{\KK}{\mathbb{K}} \newcommand{\FF}{\mathbb{F}}
\newcommand{\CC}{\mathbb{C}} \newcommand{\EE}{\mathbb{E}}
\newcommand{\XX}{X}
\newcommand{\HH}{I\hspace{-2mm}H}
\newcommand{\norm}{\Vert\hspace{-0.35mm}|}
\newcommand{\Sph}{\mathbb{S}}
\newcommand{\ganz}{\overline{\XX}}
\newcommand{\rand}{\partial\XX}
\newcommand{\prodrand}{\partial\XX_1\times\partial\XX_2} 
\newcommand{\regrand}{\partial\XX^{reg}}
\newcommand{\singrand}{\partial\XX^{sing}}
\newcommand{\Frand}{\partial^F\XX}
\newcommand{\Lim}{L_\Gamma}          
\newcommand{\Flim}{F_\Gamma}
\newcommand{\reglim}{L_\Gamma^{reg}}
\newcommand{\radlim}{L_\Gamma^{rad}}
\newcommand{\raylim}{L_\Gamma^{ray}}
\newcommand{\horinf}{\mbox{Vis}^{\infty}}
\newcommand{\horF}{\mbox{Vis}^F}
\newcommand{\Sml}{\mbox{Small}}
\newcommand{\SmlF}{\mbox{Small}^F}

\newcommand{\ifl}{\qquad\Longleftrightarrow\qquad}
\newcommand{\at}{\!\cdot\!}
\newcommand{\ging}{\gamma\in\Gamma}
\newcommand{\xo}{{o}}
\newcommand{\gamo}{{\gamma\xo}}
\newcommand{\gam}{\gamma}
\newcommand{\gax}{h}
\newcommand{\gxi}{{G\!\cdot\!\xi}}
\newcommand{\bd}{$(b,\theta)$-densit}
\newcommand{\bt}{$(b,\theta)$-densit}
\newcommand{\cd}{$(\alpha,\Gamma\at\xi)$-density}
\newcommand{\be}{\begin{eqnarray*}}
\newcommand{\ee}{\end{eqnarray*}}

\newcommand{\an}{\ \mbox{and}\ }
\newcommand{\as}{\ \mbox{as}\ }
\newcommand{\diam}{\mbox{diam}}
\newcommand{\is}{\mbox{Is}}
\newcommand{\Ax}{\mbox{Ax}}
\newcommand{\Fix}{\mbox{Fix}}
\newcommand{\Par}{F}
\newcommand{\Min}{\mbox{Fix}}
\newcommand{\rel}{\mbox{Rel}_\Gamma}
\newcommand{\vol}{\mbox{vol}}
\newcommand{\Td}{\mbox{Td}}
\newcommand{\piF}{\pi^F}

\newcommand{\for}{\ \mbox{for}\ }
\newcommand{\pr}{\mbox{pr}}
\newcommand{\sh}{\mbox{sh}}
\newcommand{\shi}{\mbox{sh}^{\infty}}
\newcommand{\rank}{\mbox{rank}}
\newcommand{\supp}{\mbox{supp}}
\newcommand{\mass}{\mbox{mass}}
\newcommand{\kernel}{\mbox{kernel}}
\newcommand{\st}{\mbox{such}\ \mbox{that}\ }
\newcommand{\Stab}{\mbox{Stab}}
\newcommand{\Root}{\Sigma}
\newcommand{\Cone}{\mbox{C}}
\newcommand{\wrt}{\mbox{with}\ \mbox{respect}\ \mbox{to}\ }
\newcommand{\where}{\ \mbox{where}\ }

\newcommand{\thet}{\widehat H}

\newcommand{\con}{{\sc Consequence}\newline}
\newcommand{\rem}{{\sc Remark}\newline}
\newcommand{\prf}{{\sl Proof}}
\newcommand{\qed}{$\hfill\Box$}

\newenvironment{rmk} {\newline{\sc Remark.\ }}{}  
\newenvironment{rmke} {{\sc Remark.\ }}{}  
\newenvironment{rmks} {{\sc Remarks.\ }}{}  
\newenvironment{nt} {{\sc Notation}}{}  

\newtheorem{satz}{\bf Theorem}

\newtheorem{df}{\sc Definition}[section]
\newtheorem{cor}[df]{\sc Corollary}
\newtheorem{thr}[df]{\bf Theorem}
\newtheorem{lem}[df]{\sc Lemma}
\newtheorem{prp}[df]{\sc Proposition}
\newtheorem{ex}{\sc Example}
\newenvironment{pros}{{\sc Properties:}}


\title{{\sc Higher order Dehn functions for 
horospheres in products of Hadamard spaces}}
\author{\sc Gabriele Link}
\date{\today}
\maketitle
\begin{abstract} Let $\XX$ 
be a product of $r$ locally compact Hadamard spaces. In this note we prove that the horospheres in $\XX$ centered at regular boundary points of $\XX$ are Lipschitz-$(r-2)$-connected. Using the filling construction by R.~Young in \cite{MR3268779} this gives sharp bounds on higher order Dehn functions for such horospheres. Moreover, if $\Gamma\subset\is(\XX)$ is a lattice acting cocompactly on $\XX$ minus a union of disjoint horoballs, we get a sharp bound on higher order Dehn functions for $\Gamma$.  We therefore deduce  that apart from the Hilbert modular groups already considered by R.~Young every irreducible $\QQ$-rank one lattice acting on a product of $r$ symmetric spaces of the noncompact type 
is undistorted up to dimension $r-1$ and has $k$-th order Dehn function asymptotic to $V^{(k+1)/k}$ for all $k\le r-2$.

\end{abstract}

\vspace{0.2cm}

\section{Introduction}

In this note we apply the filling construction introduced by R.~Young in \cite{MR3268779}  to horospheres in products of arbitrary locally compact Hadamard spaces. The only new result we need for the construction is our Proposition~\ref{slicesingbilip} which implies that the so-called "slices" in every horosphere ${\cal H}$ centered at a regular boundary point are  bilipschitz-equivalent to a Hadamard space. Using these slices one can then construct in an identical way as in \cite[Lemma~3.3]{MR3268779} a map from the $(r-1)$-simplex $\Delta^{(r-1)}\to {\cal H}$ with 
properties 
that assure Lipschitz-$(r-2)$-connectivity of ${\cal H}$. Theorem~1.3 of \cite{MR3268779} then gives the following \\[3mm]
 {\bf Theorem A}$\quad$ {\sl Any horosphere centered at a regular boundary point of a product of $r$ locally compact Hadamard spaces is undistorted up to dimension $r-1$. Moreover, for any $k\le r-2$ the $k$-th order Dehn function of such a horosphere is asymptotic to $V^{(k+1)/k}$. 
} \\[-1mm]

Moreover, if $\Gamma\subset\is(\XX)$ is a lattice acting cocompactly on $\XX$ minus 
a union of disjoint horoballs centered at regular boundary points, we 
get a sharp bound on higher order Dehn functions for $\Gamma$.  According to Section~2.4 in \cite{MR2079992}  every irreducible\break $\QQ$-rank one lattice $\Gamma$ acting on  a product $\XX$ of symmetric spaces of the noncompact type 
satisfies this assumption. We therefore get the following \\[3mm]
{\bf Theorem B}$\quad$ {\sl Let $\Gamma$ be an irreducible $\QQ$-rank one lattice acting on a product $\XX$ of $r$ symmetric spaces of the noncompact type. Then $\Gamma$ is undistorted up to dimension $r-1$. Moreover, for any $k\le r-2$ the $k$-th order Dehn function of $\,\Gamma$ is asymptotic to $V^{(k+1)/k}$. 
} \\[-1mm]

We remark that as soon as the symmetric space $\XX$ possesses at least one rank one factor $\XX_i$, then by Theorem~13.19 in~\cite{Raghunathan} 
every irreducible lattice $\Gamma\subset\is(\XX)$ has $\QQ$-rank one.



\section{Preliminaries}\label{Prelim}

The purpose of this section is to introduce terminology and notation and to summarize basic results about (products of) Hadamard spaces and their horospheres. The main references here  are \cite{MR1744486} and  \cite{MR1377265} (see also 
\cite{MR823981}).

Let $(\XX,d)$ be a metric space. A {\hl geodesic path} joining $x\in\XX$ to $y\in\XX$  is a map $\sigma$ from a closed interval $[0,l]\subset \RR$ to $\XX$ \st $\sigma(0)=x$, 
$\sigma(l)=y$ and $d(\sigma(t), \sigma(t'))=|t-t'|$ for all $t,t'\in [0,l]$.  We will denote such a geodesic path $\sigma_{x,y}$. $\XX$ is called {\hl geodesic} if any two points 
in $\XX$ can be connected by a geodesic path; if this path is unique we say that $\XX$ is {\hl uniquely geodesic}. A {\hl geodesic ray} in $\XX$ is a map $\sigma:[0,\infty)\to \XX$ such that for all $t'>t>0$ $\sigma\ein_{[t,t']}$ is a geodesic path; a  {\hl geodesic (line)} in $\XX$ is a map $\sigma:\RR\to \XX$ such that the above holds for all $t'>t$. Two geodesic rays are called {\hl asymptotic} if they are at bounded Hausdorff distance from each other.

A metric space $(\XX,d)$ is called a {\hl Hadamard space} if it is complete, geodesic and if all triangles satisfy the CAT$(0)$-inequality. This implies in particular that $\XX$ is simply connected and uniquely geodesic.  From here on we assume that  $\XX$ is a  {\hl locally compact} Hadamard space. The geometric boundary $\rand$ of $\XX$ is defined as  the set of equivalence classes of geodesic rays, where two geodesic rays are equivalent if they are asymptotic. Notice that for $x\in\XX$ and $\xi\in\rand$ arbitrary there exists a unique geodesic ray emanating 
from $x$ which belongs to the class of $\xi$; we will denote this ray by $\sigma_{x,\xi}$. Moreover, since $\XX$ is locally compact the set $\ganz:=\XX\cup\rand$ endowed with the cone topology (see \cite[Chapter II]{MR1377265}) is a compactification of $\XX$.    

Let $x, y\in \XX$, $\xi\in\rand$ and $\sigma$ a geodesic ray in the class of $\xi$. We set 
\begin{equation}\label{buseman}
 \bs_{\xi}(x, y)\,:= \lim_{s\to\infty}\big(d(x,\sigma(s))-d(y,\sigma(s))\big).
\end{equation}
This number is independent of the chosen ray $\sigma$, and the function
\begin{align*} \bs_{\xi}(\cdot , y): \quad \XX &\to  \RR,
\quad 
x \mapsto  \bs_{\xi}(x, y)\end{align*}
is called the {\hl Busemann function} centered at $\xi$ based at $y$ (see also \cite{MR1377265}, chapter~II). From the definition one immediately gets the following 
properties of the Busemann function which we will need in the sequel:
For all  $x,y,z\in\XX$ and $\xi\in\rand$ the {\hl cocycle identity}
\begin{equation}\label{cocycleid}
 \bs_{\xi}(x, z)=\bs_{\xi}(x, y)+\bs_{\xi}(y,z)
 \end{equation}
 holds, and we have 
 \begin{equation}\label{boundedness}
 |\bs_{\xi}(x, y)|\le d(x,y)\end{equation}
and $\bs_{\xi}(x, y)=d(x,y)\, $ if and only if  $y$ is a point on the geodesic ray $\sigma_{x,\xi}$.

For $\xi\in\rand$ and $x\in\XX$, the {\hl horoball} centered at $\xi$ based at $x$ is defined as the set
\[
\{ y\in\XX\mid \bs_\xi(x,y)<0\};\]
its boundary
\[{\cal H}_\xi(x):=\{ y\in\XX\mid \bs_\xi(x,y)=0\}\]
in $\XX$ is called the {\hl horosphere}  centered at $\xi$ based at $x$. 
 If $\sigma$ is a geodesic line, 
then the {\hl projection along horospheres} $p_\sigma$ is defined by
\[ p_\sigma:\XX\to\sigma(\RR),\quad x\mapsto p_\sigma(x)=\sigma\bigl(\bs_{\sigma(\infty)}(\sigma(0),x)\bigr) .\]
Notice that  $p_\sigma(x)$ is the unique intersection point of $\sigma(\RR)$ with the  horosphere ${\cal H}_{\sigma(\infty)}(x)$ centered at $\sigma(\infty)$ based at $x$.

In this note we consider the Cartesian product $\XX=\XX_1\times \XX_2\times\cdots\times\XX_r\,$ of $r$ locally compact Hadamard spaces  $(\XX_1,d_1)$, $(\XX_2,d_2),\ldots,$ $(\XX_r,d_r)$ 
endowed with the 
distance $d=\sqrt{d_1^2+d_2^2+\cdots+d_r^2}$. Notice that $(\XX,d)$ 
is again a locally compact Hadamard space.  

We denote $p_i:\XX\to \XX_i$, $i\in\{1,2,\ldots,r\}$, the natural projections. Every geodesic path $\sigma:[0,l]\to\XX\,$ can be written as a product $\sigma(t)=(\sigma_1(t\cdot \theta_1), \sigma_2(t\cdot \theta_2),\ldots, \sigma_r(t \cdot \theta_r))$, where  $\sigma_i$ are geodesic paths in $\XX_i$, $i=1,2,\ldots r$, and the $\theta_i\ge 0$ satisfy 
\[\sum_{i=1}^r \theta_i^2=1.\]  
The unit vector 
\[ \text{sl}(\sigma):= \left(\begin{array}{c} \theta_1\\ \theta_2\\\vdots\\ \theta_r\end{array}\right)\in E:=\{\theta\in  \RR^r  :\ \Vert \theta\Vert =1, \ \theta_i \ge 0\ \text{for all } i\in\{1,2,\ldots,r\}\}\]
is called the {\hl slope of $\sigma$};  a geodesic path $\sigma$ is said to be {\hl regular} if its slope does not possess a coordinate zero, i.e. if 
\[ \text{sl}(\sigma)\in E^+:=\{\theta\in  \RR^r  :\ \Vert \theta\Vert =1, \ \theta_i > 0\ \text{for all } i\in\{1,2,\ldots,r\}\},\]
and  {\hl singular} otherwise. In other words, $\sigma$ is regular if none of the projections $p_i\bigl(\sigma([0,l])\bigr)$, $i\in\{1,2,\ldots, r\}$,   is a point.

It is an easy exercise to verify that two asymptotic geodesic rays $\sigma$ and $\sigma'$  necessarily have the same slope. 
So the slope $\text{sl}(\tilde\xi)$ of a point $\tilde\xi\in\rand$ can be defined as the slope of an arbitrary geodesic ray representing $\tilde\xi$. The {\hl regular geometric boundary} $\regrand$ and  the {\hl singular geometric boundary} $\singrand$ of $\XX$ are 
defined by 
 \[ \regrand:=\{\tilde\xi\in\rand:\ \text{sl}(\tilde\xi)\in E^+\},\quad \singrand:=\rand\setminus\singrand;\]
notice that the singular boundary $\singrand$ consists of equivalence classes of geodesic rays in $\XX$ which project to a point in at least one of the factors $\XX_i$.
Moreover, given\break $\theta=(\theta_1,\theta_2,\ldots,\theta_r)\in E$ we can  define the subset
\begin{equation}\label{randtheta}
 \rand_\theta:=\{\tilde\xi\in\rand:\ \text{sl}(\tilde\xi)=\theta\}
 \end{equation}
of the geometric boundary which is 
homeomorphic to the Cartesian product of the geometric boundaries $\rand_i$ with $i\in I^+(\theta):=\{i\in \{1,2,\ldots,r\}: \theta_i>0\}$. So if $\tilde\xi\in\rand_\theta$, then for all $i\in I^+(\theta)$ the projection $\xi_i$ of $\tilde\xi$ to $\rand_i$ is well-defined.
Notice that in the particular case $\theta\in E^+$ we have $\rand_\theta\subset\regrand$, $I^+(\theta)=\{1,2,\ldots,r\}$ and hence $\rand_\theta$ is homeomorphic to 
$\rand_1\times\rand_2\times\cdots\times \rand_r$. 

The following easy lemma relates the Busemann function~(\ref{buseman}) of the product to the Busemann functions on the factors.  For a proof we refer the reader to Lemma~3.3 in \cite{1403.4858}. 
\begin{lem}\label{busprod} Let $\theta=(\theta_1,\theta_2,\ldots,\theta_r)\in E$,    $x=(x_1,x_2,\ldots,x_r)$, $y=(y_1,y_2,\ldots,y_r)\in\XX$ and $\tilde\xi\in\rand_\theta$. If $\xi_i$ denotes the projection of $\tilde\xi$ to $\rand_i$ then 
\[\bs_{\tilde\xi}(x,y)=\sum_{i\in I^+(\theta)} \theta_i \cdot \bs_{\xi_i}(x_i,y_i).\]
%
\end{lem}

\section{Horospheres and slices}

For the remainder of this note 
we fix a base point $\xo=(\xo_1,\xo_2,\ldots, \xo_r)\in X$, a point  $\tilde\xi\in\rand$ of slope $\theta\in E$ 
and abbreviate ${\cal H}:={\cal H}_{\tilde\xi}(\xo)$.

We recall the definition of  
\[ I^+(\theta)=\{i\in\{1,2,\ldots,r\}:\ \theta_i>0\} \]
and remark that for all $i\in I^+(\theta)$ the projection $\xi_i\,$ of $\tilde\xi\,$ to $\rand_i$ is well-defined. More\-over, 
  according to Lemma~\ref{busprod} we have
\[ {\cal H}=\{y=(y_1,y_2,\ldots,y_r)\in X: \ \bs_{\tilde\xi}(\xo,y)=\sum_{i\in I^+(\theta)} \theta_i \cdot \bs_{\xi_i}(\xo_i,y_i)=0\}.\]
Let $q\in \{1,\ldots, r\}$ denote the cardinality of the set $I^+(\theta)$. Notice that if $\tilde\xi\in\singrand$ then $q\le r-1$. 

As in \cite{MR3268779} we consider the following important subsets of the horosphere ${\cal H}$: 
For  $i\in I^+(\theta)$ we let $s_i\subseteq \XX_i\,$ be either all of $\XX_i\,$ or the image of a geodesic line $\sigma_i\,$ in $X_i$ with $\sigma_i(\infty)=\xi_i$. For $i\notin I^+(\theta)$ we set $s_i=\{\xo_i\}$. If precisely $k$ of the sets $s_i$ are all of $\XX_i$, then the set 
\[ (s_1\times s_2\times\cdots\times s_r)\cap {\cal H} \]
is called a {\hd $k$-slice}; notice that necessarily $k\in\{0,1,\ldots, q-1\}$.  For a $k$-slice\break $S= (s_1\times s_2\times\cdots\times s_r)\cap {\cal H} $ we denote  
\begin{equation}\label{indfullfactors} I(S):=\{ i\in\{1,2,\ldots, r\} :\ s_i=\XX_i\}\end{equation}
the index set of the "full factors"  which by definition has cardinality $k$.

We now state the key proposition for horospheres ${\cal H}$ centered at a boundary point $\tilde\xi\in \rand_\theta$; recall that $q=\# I^+(\theta)$:
\begin{prp}\label{slicesingbilip} If $k\in\{0,1,\ldots, q-1\}$ and $S=(s_1\times s_2\times\cdots\times s_r)\cap {\cal H}\subset {\cal H}$ is a $k$-slice, then $S$ together with the length metric $d_S$ induced from the metric $d$ on $\XX$ is bilipschitz-equivalent to the Cartesian product of the factors $\XX_i$ with $i\in I(S)$ times $\RR^{q-k-1}$ (endowed with the product metric). 
\end{prp}

\prf.\ \ For simplicity of notation we first assume that 
$I^+(\theta)= \{1,2,\ldots,q\}$ and hence
\[ {\cal H}=\{ y=(y_1,y_2,\ldots,y_r)\in X: \ \bs_{\tilde\xi}(\xo,y)=\sum_{i=1}^q \theta_i \cdot \bs_{\xi_i}(\xo_i,y_i)=0\}.\]
So for a $k$-slice the last $r-q\,$ of the $s_i$ are always equal to $\{\xo_i\}$, and we moreover assume that the first $k$ of the $s_i$ are equal to $\XX_i$.  For $l\in\{k+1,\ldots, q\}$ we let $\hat \sigma_l\subset\XX_l$ be a geodesic line with $\hat\sigma_l(\RR)=s_l$.  Reparametrizing $\hat\sigma_l$ if necessary we may further require that
$\bs_{\xi_l}(\xo_l, \hat\sigma_l(0))=0$. We prove that the map
\begin{align*} \text{pr}&:\quad S\to \XX_1\times\cdots\times\XX_k\times \RR^{q-k-1},\\
&\quad (x_1,x_2,\ldots,x_k,x_{k+1},\ldots, x_q,\xo_{q+1},\ldots \xo_{r})\\
& \qquad\qquad \mapsto \bigl(x_1,\ldots, x_k,\bs_{\xi_{k+1}}(\xo_{k+1},x_{k+1}),\ldots, \bs_{\xi_{q-1}}(\xo_{q-1},x_{q-1})\bigr) \end{align*}
is the desired bilipschitz-equivalence.

We first show that $\text{pr}$ is bijective: \\
For that we let $(x_1,\ldots,x_k, t_{k+1},\ldots, t_{q-1})\in  \XX_1\times\cdots\times\XX_k\times \RR^{q-k-1}$ be arbitrary. In order to construct the preimage by $\text{pr}$ we first set $x_l= \hat\sigma_l(t_l)\,$ for $l\in\{k+1,\ldots, q-1\}$; in particular we have
\[ \bs_{\xi_l}(\xo_l, x_l)= \bs_{\xi_l}\bigl(\xo_l, \hat\sigma_l(t_l)\bigr)= \underbrace{\bs_{\xi_l}\bigl(\xo_l, \hat\sigma_l(0)\bigr)}_{=0} + \bs_{\xi_l}\bigl(\hat\sigma_l(0), \hat\sigma_l(t_l)\bigr)=t_l\]
since $\hat\sigma_l(\infty)=\xi_l$. 
Finally we set 
\[ t_q:=-\frac1{\theta_q} \Bigl( \sum_{i=1}^{q-1}\theta_i \bs_{\xi_i}(\xo_i,x_i)\Bigr)\]
 and $x_q:=  \hat\sigma_l(t_q)$. 

 We next prove that $\text{pr}$ is bilipschitz. For that we will denote the product metric on $ \XX_1\times\cdots\times\XX_k\times \RR^{q-k-1}$ by $\overline d$. Let  $x=(x_1,x_2,\ldots,x_r)$, $y=(y_1,y_2,\ldots,y_r)\in S$  be arbitrary. Since $x_i=y_i=\xo_i\,$ for $i\in \{q+1,\ldots r\}$ we clearly have
 \begin{align*}
 d_S(x,y)^2\ge d_X(x,y)^2&=\sum_{i=1}^r d_i(x_i,y_i)^2=\sum_{i=1}^q d_i(x_i,y_i)^2\\
 &\ge \sum_{i=1}^{q-1} d_i(x_i,y_i)^2=\overline{d}\bigl(\pr(x),\pr(y)\bigr)^2,
 \end{align*}
 because for $k+1\le i\le q-1$ the points $x_i$ and $y_i$ lie on the geodesic line $\hat\sigma_i$.
 
 In order to prove the converse estimate, we will construct a path $c$ in $S$ joining $x$ and $y$ with length bounded by a constant times $\overline{d}\bigl(\pr(x),\pr(y)\bigr)$. We abbreviate
 \[  \overline{d}=\overline{d}\bigl(\pr(x),\pr(y)\bigr)\qquad\text{and}\quad  d_i=d_i(x_i,y_i)\qquad\text{for }\  1\le i\le q.\]
  For $1\le i\le q-1$  we consider the geodesic segment $\sigma_i:[0,d_i]\to\XX_i$ joining $x_i$ and $y_i$. Then there exists a unique curve $\gamma: [0,\overline d]\to\RR\,$ \st the path
 \[ c:[0,\overline d]\to\XX,\quad t\mapsto \bigl(\sigma_1(t\theta_1),\sigma_2(t\theta_2),\ldots,\sigma_{q-1}(t\theta_{q-1}), \hat\sigma_q(\gamma(t)),\xo_{q+1},\ldots,\xo_{r}\bigr)\]
 is contained in the $k$-slice $S$. Indeed, the curve $\gamma$ is determined by the condition
 \begin{equation}\label{condthetaq} \theta_q\bs_{\xi_q}\bigl(\xo_q,\hat\sigma_q(\gamma(t))\bigr)=-\sum_{i=1}^{q-1} \theta_i\bs_{\xi_i}(\xo_i,\sigma_i(t\theta_i))\quad\text{for all }\ t\in [0,\overline d].\end{equation}
Using the properties 
of the Busemann functions as well as the Cauchy-Schwarz inequality, we estimate for $t,t'\in [0,\overline d]$ 
 \begin{align*}
| \gamma(t')-\gamma(t)|&
 =\big| \bs_{\xi_q}(\hat\sigma_q(\gamma(t)),\hat\sigma_q(\gamma(t'))\big|
\stackrel{(\ref{cocycleid})}{=}\big| \bs_{\xi_q}(\hat\sigma_q(\gamma(t)),\xo_q)+\bs_{\xi_q}(\xo_q,\hat\sigma_q(\gamma(t'))\big|\\
 &\stackrel{(\ref{condthetaq})}{=} \frac1{\theta_q}\Big| \sum_{i=1}^{q-1} \theta_i\bigl( \bs_{\xi_i}(\sigma_i(t\theta_i),\xo_i)+\bs_{\xi_i}(\xo_i,\sigma_i(t'\theta_i))\bigr)\Big|\\
 &\stackrel{(\ref{cocycleid})}{\le} \frac1{\theta_q}\sum_{i=1}^{q-1} \theta_i \big|\bs_{\xi_i}(\sigma_i(t\theta_i),\sigma_i(t'\theta_i))\big|\stackrel{(\ref{boundedness})}{\le} \frac1{\theta_q} \sum_{i=1}^{q-1} \theta_i d_i\bigl(\sigma_i(t\theta_i),\sigma_i(t'\theta_i)\bigr)\\
 &\le \frac1{\theta_q}\sqrt{ \sum_{i=1}^{q-1}\theta_i^2}\sqrt{ \sum_{i=1}^{q-1}  d_i\bigl(\sigma_i(t\theta_i),\sigma_i(t'\theta_i)\bigr)^2}= \frac{\sqrt{1-\theta_q^2}}{\theta_q}\sqrt{\sum_{i=1}^{q-1} |t-t'|^2\theta_i^2}\\
 &= |t-t'| \frac{1-\theta_q^2}{\theta_q}\le \frac1{\theta_q} |t-t'| .
  \end{align*}
  We conclude that the length of the curve $\hat\sigma_q\circ \gamma\,$ in $\XX_q$  (which is equal to the length of the curve $\gamma$ in $\RR$) is smaller than or equal to $\overline d/\theta_q$. So for the length $L(c)$ of the curve $c$ in $\XX$ we get
  \[ L(c)\le \sqrt{\sum_{i=1}^{q-1} d_i^2+\Bigl(\frac{\overline{d}}{\theta_q} \Bigr)^2}= \overline d\, \sqrt{1+\frac1{\theta_q^2}},\] 
  and hence \[ d_S(x,y)\le \sqrt{1+\frac1{\theta_q^2}}\cdot\overline d(\pr(x),\pr(y).\]
 
Finally, for $I^+(\theta)$ arbitrary of cardinality $q$  and any $k$-slice $S\subset{\cal H}$ with index set $I(S)\subset\{1,2,\ldots, r\}$ of the "full factors" in $S$  
we denote by $\pr$ the bijection from $S$ to the Cartesian product of the factors $\XX_i$ with $i\in I(S)$ times $\RR^{q-k-1}$, and $\overline{d}$  its product metric. Then for all $x,y\in S$ we have
 \[ \overline{d}(\pr(x),\pr(y))\le d_S(x,y)\le  \sqrt{1+\max\Big\{\frac1{\theta_i^2}:\, i\in I^+(\theta) \Big\}}\,\overline{d}(\pr(x),\pr(y)),\]
 which proves the claim.\qed\\

\section{Statement of results}

In order to deduce Theorems~A and B from the introduction we use Lemma~3.2 in \cite{MR3268779}, which is due to Gromov and provides an equivalent formulation of Lipschitz connectivity. More precisely, if  ${\cal H}$ is a horosphere centered at a boundary point $\tilde\xi\in\rand_\theta$ with $q=\# I^+(\theta)$, we construct a map from the $(q-1)$-simplex $\Delta^{(q-1)}$ to the horosphere ${\cal H}$ with the 
properties (1) and (2) of Lemma~3.2 in \cite{MR3268779}. This can be done word by word as in Section~3 of the article by R.~Young, since our Proposition~\ref{slicesingbilip} ensures that for $k\in\{0,1,\ldots q-1\}$ every $k$-slice of the horosphere ${\cal H}$  is bilipschitz-equivalent to a Hadamard space. This proves that 
 ${\cal H}$ is Lipschitz-$(q-2)$-connected, hence applying Theorem~1.3 and Corollary~1.4 of \cite{MR3268779} we get the following stronger version of Theorem~A:
\begin{thr}\label{eins}
Let $\XX$ be a product of $r$ locally compact Hadamard spaces and $\tilde \xi\in\rand$. If $\tilde \xi$ is represented by a geodesic ray which does not project to a point in precisely $q$ factors, then every horosphere ${\cal H}$ centered at $\tilde \xi$ is undistorted up to dimension $q-1$. Moreover, for any $k\le q-2$ the $k$-th order Dehn function of such a horosphere is asymptotic to $V^{(k+1)/k}$.
\end{thr}



We want to remark that for horospheres centered at a singular boundary point $\tilde\xi\in\rand_\theta$, the   $k$-slices could be optionally defined by setting $s_i=\XX_i$ (instead of $s_i=\{\xo_i\}$) for $i\notin I^+(\theta)$. Mimicking our proof of Proposition~\ref{slicesingbilip} one can show that every such $k$-slice $S$ is bilipschitz-equivalent to the Cartesian product of the factors $\XX_i$ with $i\in I(S)$ (as defined by~(\ref{indfullfactors})) 
times $\RR^{q-k-1}$ via the map pr which is the identity on the product of the "full factors" $\XX_i$ with $i\in I(S)$ as in the proof of Proposition~\ref{slicesingbilip}. However, now we always have $\# I(S)\ge  \# I^+(\theta)= r-q$, so the minimal dimension of a $k$-slice would be $r-q\ge 1$; but since in Young's construction one needs to start with $0$-slices, this would not yield a statement about Lipschitz-connectivity.

We finally turn to the proof of Theorem~B. Since undistortedness and (higher order) Dehn functions are quasi-isometry invariants, Theorem~\ref{eins} implies \begin{thr}
Let $\XX$ be a product of $r$ locally compact Hadamard spaces and $\Gamma\subset\is(\XX)$  a lattice acting cocompactly on $\XX$ minus a union of disjoint horoballs. If all the centers of these disjoint horoballs are represented by a geodesic ray which does not project to a point in at least $q$ factors, then $\Gamma$ is undistorted up to dimension $q-1$. Moreover, for any $k\le q-2$ the $k$-th order Dehn function of $\,\Gamma$ is asymptotic to $V^{(k+1)/k}$.
\end{thr}

Obviously, if all of the disjoint horoballs are centered at regular boundary points then the conclusion of the above theorem holds with $q=r$. Moreover, by Prasad \cite{MR0385005} and Raghunathan \cite{Raghunathan} every irreducible lattice of $\QQ$-rank one in a semi-simple Lie group acts cocompactly on the associated symmetric space $\XX$ minus a finite disjoint union of horoballs; Proposition~2.4.3 in \cite{MR2079992} further states that the centers of these horoballs belong to the regular boundary. So Theorem~B follows.

\bibliography{Bibliographie} 

\vspace{1.5cm}
\noindent Gabriele Link\\
Institut f\"ur Algebra und Geometrie\\
Karlsruher Institut f\"ur Technologie (KIT)\\
Englerstr. 2, Geb\"aude 20.30  \\[1mm]
D-76131 Karlsruhe\\
e-mail:\ gabriele.link@kit.edu

\end{document}